\documentclass[oneside]{article}
\usepackage{amsmath}
\usepackage{amssymb}
\usepackage[mathscr]{euscript}
\usepackage{hyperref}
\usepackage{tikz}
\date{June 26, 2012}

\newtheorem{thm}{Theorem}[section]

\newtheorem{lemma}[thm]{Lemma}
\newtheorem{theorem}[thm]{Theorem}

\newtheorem{remark}[thm]{Remark}

\newtheorem{definition}[thm]{Definition}

\newtheorem{conjecture}[thm]{Conjecture}

\newenvironment{proof}{{\bf Proof:}}{\hfill$\square$\vskip.5cm}

\newcommand{\R}{\mathbb{R}}
\newcommand{\N}{\mathbb{N}}

\newcommand{\E}{\mathbf{E}}

\newcommand{\neswarrow}{\begin{tikzpicture}\draw[<->](0,0)--(0.15,0.15);\end{tikzpicture}}
\newcommand{\ind}{\hspace{-2pt}\begin{minipage}{15pt}\begin{tikzpicture}\draw[](0,0) node[rotate=90] {$\models$};\end{tikzpicture}\end{minipage}}

\title{About Thinning Invariant Partition Structures\footnote{This is a revision. The original version was written in May 2011, and we have left the affiliations
as they were then.}}

\author{Shannon Starr$^1$, Brigitta Vermesi$^2$ and Ang Wei$^1$\\[10pt]
$^1$ Department of Mathematics\\ University of Rochester\\ Rochester, NY 14627\\[5pt]
$^2$ University of Washington\\ Department of Mathematics\\ Seattle, WA 98195}

\begin{document}

\maketitle

\begin{abstract}
\setcounter{section}{0}
Bernoulli-$p$ thinning has been well-studied for point processes.
Here we consider three other cases: (1) sequences $(X_1,X_2,\dots)$;
(2) gaps of such sequences $(X_{n+1}-X_1)_{n\in\mathbb{N}}$;
(3) partition structures.
For the first case we characterize the distributions which are simultaneously
invariant under Bernoulli-$p$ thinning for all $p \in (0,1]$.
Based on this, we make conjectures for the latter two cases, and provide
a potential approach for proof.
We explain the relation to spin glasses, which is complementary to 
important previous work of Aizenman and Ruzmaikina, Arguin,
and Shkolnikov.
\end{abstract}

\section{Introduction and main conjecture} 

Motivated by the theory of mean field spin glasses, we consider thinning for random partition
structures.
In spin glasses a configuration consists of $N$ spins $(\sigma_1,\dots,\sigma_N)$ in $\{+1,-1\}^N$
and has energy proportional to $N$.
One wishes to take $N \to \infty$ to obtain the thermodynamic limit.
But
each time $N$ increases by 1, the number of configurations doubles.
Additionally there is a shift in the energy of each configuration at least of order 1, because the energy also scales as $N$.

Both the energy shift and the entropy shift must be considered in order to understand equilibrium
states from the ``cavity'' perspective.
For the entropy part of the dynamics,
it is more intrinsic to imagine taking spins away one at a time.
This leads to the notion of thinning.
Thinning has been much studied in the context of random point processes,
most notably by Matthes, Kerstan and Mecke \cite{MatthesKerstanMecke}
and Kallenberg \cite{Kallenberg}.

For uncorrelated spin glasses, such as Derrida's Random Energy Model (REM) \cite{Derrida}, the natural objects
are random partitions structures.
These are random probability measures, ignoring the underlying structure of the points of the sample space.
This structure is unimportant for the REM
because of the IID nature of the energies.

In an important paper, Ruzmaikina and Aizenman considered competing particle systems,
from the perspective of energy gaps,
motivated by the energy shift aspect of the REM \cite{RuzmaikinaAizenman}.
This was followed by alternative perspectives by Arguin \cite{Arguin}
who considered random partitions structures,
and Shkolnikov \cite{Shkolnikov} who considered a more general family of energy
shifts allowing for lattice type.

In the present paper we consider the entropy shift from the perspective of thinning for random partition structures.
We describe this setting next.
Our interest is to characterize the set of random partition structures which are thinning invariant.
In particular, we will state a conjecture and support it with evidence.


\smallskip

{\bf Notation.} We write $X \sim \mu$ to mean that $X$ is a random variable with distribution $\mu$.
We write $X \ind Y$ to indicate that $X$ and $Y$ are independent.
Given a random variable $Z$ which has been constructed in some way, we will write $\mathscr{L}(Z)$
for the marginal distribution (or law).

\subsection{Bernoulli-$p$ thinning for partition structures}
Let $\Delta$ denote the set of all $\xi = (\xi_1,\xi_2,\dots)$ satisfying
$\xi_1\geq \xi_2\geq \dots \geq 0$ and $\xi_1 + \xi_2 + \dots \leq 1$.
These are the partition structures.
With the product topology, $\Delta$ is a compact and metrizable space.

The continuous component of the random measure is $\xi_0 := 1 - (\xi_1 + \xi_2 + \dots)$.
This is sometimes also called ``dust.''
We define $\Delta_{\infty}$ to be the subset of those $\xi \in \Delta$ such that one or both of the conditions is satisfied
\begin{itemize}
\item
$\xi_0>0$, or
\item $|\{n\, :\, \xi_n>0\}|=\infty$.
\end{itemize}
We call these infinite partition structures.

Let $\mathcal{M}(\Delta)$ denote the set of Borel probability measures on $\Delta$.
These are the random partition structures.

\begin{definition}
Given a complete separable metric space $\mathscr{X}$, we always let $\mathcal{M}(\mathscr{X})$
denote the set of Borel probability measures on $\mathscr{X}$.
\end{definition}

Let $\mathcal{M}(\Delta_{\infty})$ denote the subset of those $\mu \in \mathcal{M}(\Delta)$
satisfying $\mu(\Delta_{\infty})=1$.
These are the random infinite partition structures.

For each $p \in (0,1]$, the thinning map is $\theta_p : \mathcal{M}(\Delta_{\infty}) \to \mathcal{M}(\Delta_{\infty})$,
defined as follows.
Given $\mu \in \mathcal{M}(\Delta_{\infty})$, let $\xi=(\xi_1,\xi_2,\dots)$ be a random element of $\Delta_{\infty}$,
$$
\xi \sim \mu\, .
$$
Let $B = (B_1,B_2,\dots)$ be independent Bernoulli-$p$ random variables, 
$$
B \ind \xi\, .
$$
Let $Z' = p \xi_0 + B_1 \xi_1 + B_2 \xi_2 + \dots$.
Let $K_1 = \min\{k : B_k=1\}$ and inductively define $K_{n+1} = \min\{k>K_n\, :\, B_k=1\}$,
for $n\in \N$. Note that, a.s., such a sequence exists and is unique.
Let us define $\zeta  = (\zeta_1,\zeta_2,\dots)$ such that
$$
\zeta_n\,  =\, \xi_{K_n}/Z' \ \text{ for $n \in \N$,}
$$
so that  $\zeta$ is an element of $\Delta_{\infty}$.
Then with all this
$$
\theta_p(\mu)\, \stackrel{\text{def}}{:=}\, \mathscr{L}(\zeta)\, .
$$

We say that $\mu \in \mathcal{M}(\Delta_{\infty})$
is thinning invariant if $\theta_p(\mu) = \mu$ for all $p \in (0,1]$.
Our main goal is to characterize the class of all thinning invariant partition structures,
which we will do as a conjecture.

\subsection{Poisson-Kingman partition structures}

An excellent reference for Poisson-Kingman partition structures is the review by Pitman \cite{Pitman}.
These are also related to Lambda-coalescents as in \cite{PitmanLambda}.
For a detailed review we recommend  \cite{Berestycki}.

Let $\mathcal{L}$ denote the set of all pairs $(\Lambda,v)$ 
where $v$ is a nonnegative number
and $\Lambda$ is a Borel measure on $(0,\infty)$ satisfying
$$
\int \min\{x,1\}\, d\Lambda(x)<\infty\, .
$$
We let $\mathcal{L}_{\infty}$ denote the subset consisting of all pairs $(\Lambda,v)\in \mathcal{L}$
such that one or both of the following conditions is satisfied
\begin{itemize}
\item $v>0$, or 
\item $\Lambda((0,\infty))=\infty$.
\end{itemize}
Given $(\Lambda,v) \in \mathcal{L}_{\infty}$, the Poisson-Kingman random partition 
structure $\nu_{(\Lambda,v)} \in \mathcal{M}(\Delta_{\infty})$ is defined as follows.

\begin{definition}
\label{def:MGF}
Given a complete separable metric space
$\mathscr{X}$,
let $\mathcal{M}_+(\mathscr{X})$ denote the set of all (nonnegative) Borel measures,
not necessarily normalized or finite,
and let $\mathcal{B}_+(\mathscr{X})$ denote the set of all nonnegative Borel functions.
Given $\rho \in \mathcal{M}_+(\mathscr{X})$ let $\operatorname{PPP}(\mathscr{X},\rho) \in \mathcal{M}(\mathcal{M}_+(\mathscr{X}))$ denote the
Poisson point process on $\mathscr{X}$ with intensity measure $\rho$, which 
is the distribution of a random Borel measure $\Xi \in \mathcal{M}_+(\mathscr{X})$, 
specified by its moment generating functional:
$$
\Xi\, \sim\, \operatorname{PPP}(\mathscr{X},\rho)\quad \Leftrightarrow\quad
\forall f \in \mathcal{B}_+(\mathscr{X})\, ,\ 
\E[e^{-\int_{\mathscr{X}} f\, d\Xi}]\, =\, \exp\left(-\int_{\mathscr{X}} (1-e^{-f})\, d\rho\right)\, .
$$
See Chapter 7 of Daley and Vere-Jones for more details \cite{DaleyVere-Jones}
\end{definition}

Now let $\Xi$ be a random measure on $(0,\infty)$,
$$
\Xi\, \sim\, \operatorname{PPP}((0,\infty),\Lambda)\, .
$$
We define the random variable $Z = v + \Xi((0,\infty))$.
Almost surely there exists a sequence
$\xi = (\xi_1,\xi_2,\dots)$ with $\xi_1 \geq \xi_2\geq \dots \geq 0$
and such that 
$$
\Xi\, =\, \sum_{n=1}^{\infty} \delta_{\xi_n}\, .
$$
We define $\widehat{\xi}_n = \xi_n / Z$ for each $n \in \N$,
and let $\widehat{\xi} = (\widehat{\xi}_1,\widehat{\xi}_2,\dots)$.
Then 
$$
\nu_{(\Lambda,v)}\, \stackrel{\text{def}}{:=}\, \mathscr{L}(\widehat{\xi})\, .
$$

\subsection{Conjecture for partition structures}
\label{sec:Conjecture}

For each $x \in \R$, define the shift $\sigma_x : \R \to \R$ as $\sigma_x(y) = x+y$ for each $y \in \R$.
Let $\Sigma_x : \mathcal{M}_+(\R) \to \mathcal{M}_+(\R)$ denote $\Sigma_x(\rho) = \rho \circ \sigma_x^{-1}$.
Let us define 
$$
\mathcal{M}_{\text{st}}(\mathcal{M}_+(\R))\, \stackrel{\text{def}}{:=}\,
\{Q \in \mathcal{M}(\mathcal{M}_+(\R))\, :\, \forall x \in \R\, ,\
Q \circ \Sigma_x^{-1}\, =\, Q\}\, .
$$
In other words, $\mathcal{M}_{\text{st}}(\mathcal{M}_+(\R))$ is the set of stationary random measures.

Let $|\cdot|$ denote Lebesgue measure on $\R$.
Let $\mathcal{M}_{+,1}(\R)$ denote the subset of all those $\rho \in \mathcal{M}_+(\R)$ satisfying
$$
\lim_{t \to \infty} t^{-1} \int_{0}^{t} \rho(\sigma_x^{-1}(A))\, dx\, =\, 
\lim_{t \to \infty} t^{-1} \int_{-t}^{0} \rho(\sigma_x^{-1}(A))\, dx\, =\, |A|\, 
$$
for every Borel set $A \subseteq \R$.

Given $\rho \in \mathcal{M}_{+,1}(\R)$ and given $m>0$, let us first define $\rho^{(m)} \in \mathcal{M}_+(\R)$ as
$$
d\rho^{(m)}(x)\, =\, m e^{mx}\, d\rho(x)\, ,
$$
and then let us define
$\Lambda_{\rho,m} \in \mathcal{M}_+((0,\infty))$ by
$$
\Lambda_{\rho,m}\, =\, \rho^{(m)} \circ e_{-}^{-1}\, ,
$$
where $e_{-} : \R \to (0,\infty)$ is the mapping $e_{-}(x) = e^{-x}$.

\begin{remark}
Consider the special case $\rho = |\cdot|$. Then $\Lambda_{\rho,m}(dx) = m x^{-m-1}\, dx$ for each $m \in (0,1)$.
Then the Poisson-Kingman partition structure $\nu_{(\Lambda_{\rho,m},0)}$ is called the Poisson-Dirichlet 
random partition structure $\operatorname{PD}(m,0) \in \mathcal{M}(\Delta_{\infty})$
by Pitman and Yor \cite{PitmanYor}.
This is important in spin glass theory as in \cite{RuzmaikinaAizenman,Arguin}.
\end{remark}

Let $\mathscr{K}$ denote the disjoint union of $\mathcal{M}_{+,1}(\R) \times (0,1)$ and $\{1\}$.
Let $\mathcal{M}_{\text{st}}(\mathscr{K})$ denote the set  of Borel probability measures on $\mathscr{K}$
such that for each Borel set $A\subseteq \mathcal{M}_{+,1}(\R) \times (0,1)$ and each $x \in \R$,
$$
Q(\{(\rho,m)\, :\, (\Sigma_x(\rho),m) \in A\})\, =\, Q(A)\, .
$$
Given $Q \in \mathcal{M}_{\text{st}}(\mathscr{K})$, we define $\mathcal{N}_Q \in \mathcal{M}(\Delta_{\infty})$ to be a special
example of a ``Cox-Kingman''
partition structure, defined such that
for any Borel set $A \subseteq \Delta_{\infty}$,
$$
\mathcal{N}_Q(A)\, =\, Q(\{1\}) \nu_{(0,1)}(A) + \int_{\mathcal{M}_{+,1}(\R)\times (0,1)} \nu_{(\Lambda_{\rho,m},0)}(A)\, dQ(\rho,m)\, .
$$

\begin{conjecture}
\label{conj:main}
(i) For any pair $Q,Q' \in \mathcal{M}_{\text{st}}(\mathscr{K})$, if $\mathcal{N}_Q = \mathcal{N}_{Q'}$ then $Q=Q'$.\\[5pt]
(ii) The set of thinning invariant $\mu \in \mathcal{M}(\Delta_{\infty})$ is precisely $\{\mathcal{N}_Q\, :\, Q \in \mathcal{M}_{\text{st}}(\mathscr{K})\}$.
\end{conjecture}


\subsection{Outline for the rest of the paper}
In the next section we are going to state the analogous problem of thinning for sequences. For this problem we have a rigorous characterization.
This is also an interesting question on its own for reasons we will explain.
In Section \ref{sec:gap} we will state the analogous problem for ``gaps.'' In a way this generalizes the set of thinning invariant random partition
structures. We also obtain a conjecture, which we support by an argument parallel to the theorem for thinning invariant sequences.
But there is one missing step of our argument which is establishing tightness of a certain sequence of measures.
Finally, in Section \ref{sec:return} we return to the problem of thinning invariant random partition structures.
In particular, we discuss the connection to Derrida's REM more explicitly than in the introduction, showing the two
complementary pieces of the cavity step involving the energy shift that Aizenman and Ruzmaikina already solved
and the entropy shift corresponding to thinning.

\section{Thinning-invariant sequences}

Let $\mathcal{U}$ be a compact, metrizable space.
Let $\mathcal{U}^{\N}$ denote the set of all sequences $u = (u_1,u_2,\dots)$ with
$u_1,u_2,\dots \in \mathcal{U}$.
Given $\mu \in \mathcal{M}(\mathcal{U}^{\N})$ and $p \in (0,1]$ we consider
a new measure $\Theta_p(\mu)$ in $\mathcal{M}(\mathcal{U}^{\N})$, defined as follows.
Let $U = (U_1,U_2,\dots)$ be a random element of $\mathcal{U}^{\N}$, 
$$
U \sim \mu\, .
$$
Additionally, let $B = (B_1,B_2,\dots)$ be i.i.d., Bernoulli-$p$ random variables,
$$
B \ind U\, .
$$
Let $K_1<K_2<\dots$ be the sequence uniquely defined, a.s., by the condition 
$$
\{K_1,K_2,\dots\}\, =\, \{k \in \N\, :\, B_k=1\}\, .
$$
We define $V = (V_1,V_2,\dots)$ as $V_n = U_{K_n}$ for $n \in \N$.
With all this, we define
$$
\Theta_p(\mu)\, \stackrel{\text{def}}{:=}\, \mathscr{L}(V)\, .
$$
We say that $\mu \in \mathcal{M}(\mathcal{U}^{\N})$ 
is thinning invariant if $\Theta_p(\mu)=\mu$ for all $p \in (0,1]$.

\subsection{Motivation}

After we initially stated our conjecture for thinning invariant random partition structures Dmitry Panchenko suggested characterizing thinning invariant sequences.
Panchenko also told us of the connection to exchangeability via ``spreadability.''
Ryll-Nardzewski introduced the notion of spreadability \cite{Ryll-Nardzewski}.
A sequence is spreadable if, for every non-random subsequence $k_1<k_2<\dots$,
$$
(X_{k_1},X_{k_2},\dots)\, \stackrel{\mathcal{D}}{=}\, (X_1,X_2,\dots)\, .
$$
It is easy to see that this is equivalent to exchangeability.
Our main theorem below will show that the notion of thinning invariance is strictly weaker than spreadability.

\subsection{The Poisson construction for sequences}

Let $\mathcal{C}_{\neswarrow}(\R)$ denote the set of all strictly increasing homeomorphisms $\phi : \R \to \R$.
Let $\mathcal{M}_{\neswarrow}(\R)$ denote the subset of all $\alpha \in \mathcal{M}_+(\R)$ such that the mapping
$$
F_{\alpha}(x)\, =\, \alpha((-\infty,x))
$$ 
is a homeomorphism of $\R$ onto $(0,\infty)$.
There is a right action of $\mathcal{C}_{\neswarrow}(\R)$ on $\mathcal{M}_{\neswarrow}(\R)$:
$\alpha \mapsto \alpha \circ \phi^{-1}$.

Let $\mathcal{M}_{\neswarrow}(\R \times \mathcal{U})$ denote the set of all Borel measures $\alpha$
on $\R \times \mathcal{U}$ such that the marginal $\alpha(\cdot \times \mathcal{U})$ is in $\mathcal{M}_{\neswarrow}(\R)$.
Given $\alpha \in \mathcal{M}_{\neswarrow}(\R\times \mathcal{U})$ we define $\widetilde{\nu}_{\alpha} \in \mathcal{M}(\mathcal{U}^{\N})$
as follows.
Let $\Xi$ be a random measure on $\R \times \mathcal{U}$,
$$
\Xi\, \sim\, \operatorname{PPP}(\R \times \mathcal{U},\rho)\, .
$$ 
Almost surely there is a unique random sequence $(X_1,U_1), (X_2,U_2), \dots \in \R \times \mathcal{U}$
such that $X_1<X_2<\dots$ and 
$$
\Xi\, =\, \sum_{n=1}^{\infty} \delta_{(X_n,U_n)}\, .
$$
Let $U=(U_1,U_2,\dots)$, and define
$$
\widetilde{\nu}_{\alpha}\, \stackrel{\text{def}}{:=}\, \mathscr{L}(U)\, .
$$

\subsubsection{Symmetry}
There is still a right action of $\mathcal{C}_{\neswarrow}(\R)$ on $\mathcal{M}_{\neswarrow}(\R \times \mathcal{U})$: 
$\alpha \mapsto \alpha \circ (\phi \times \operatorname{id}_{\mathcal{U}})^{-1}$ where $\operatorname{id}_{\mathcal{U}}$ is the identity mapping on $\mathcal{U}$.
Inspection shows that 
$$
\widetilde{\nu}_{\rho \circ (\phi \times \operatorname{id}_{\mathcal{U}})^{-1}}\, 
=\, \widetilde{\nu}_{\rho}\,  .
$$
There are distinguished representatives for $\mathcal{M}_{\neswarrow}(\R \times \mathcal{U})/\mathcal{C}_{\neswarrow}(\R)$.

Let $\mathcal{M}_{\operatorname{Leb}}(\R\times \mathcal{U})$ denote the set of all Borel measures $\gamma$ on $\R \times \mathcal{U}$
such that the marginal $\gamma(\cdot \times \mathcal{U})$ is $|\cdot|$.
Given $\gamma \in \mathcal{M}_{\operatorname{Leb}}(\R\times \mathcal{U})$, we define the new Borel measure $\gamma^{(1)} \in \mathcal{M}_{\neswarrow}(\R\times \mathcal{U})$
by taking 
\begin{equation}
\label{eq:1prime}
\gamma^{(1)}(dx\times du)\, =\, e^{x} \gamma(dx\times du)\, .
\end{equation}

\subsection{Cox construction and characterization}

Recall that for each $x \in \R$ we defined $\sigma_x \in \mathcal{C}_{\neswarrow}(\R)$ as $\sigma_x(\cdot) = \cdot+x$.
Define
$$
\widetilde{\Sigma}_x : \mathcal{M}_{\text{Leb}}(\R \times \mathcal{U}) \to \mathcal{M}_{\text{Leb}}(\R \times \mathcal{U})\, ,\qquad
\widetilde{\Sigma}_x(\gamma)\, =\, \gamma \circ (\sigma_x \times \operatorname{id}_{\mathcal{U}})^{-1}\, .
$$
Then we denote the stationary distributions for random measures
$$
\mathcal{M}_{\text{st}}(\mathcal{M}_{\operatorname{Leb}}(\R\times \mathcal{U}))\,
=\, \{\widetilde{Q} \in \mathcal{M}(\mathcal{M}_{\text{Leb}}(\R \times \mathcal{U}))\, :\,
\forall x \in \R\, ,\ \widetilde{Q} \circ \widetilde{\Sigma}_x^{-1} = \widetilde{Q}\}\, . 
$$
Given $\widetilde{Q} \in \mathcal{M}_{\text{st}}(\mathcal{M}_{\operatorname{Leb}}(\R\times \mathcal{U}))$, we define $\widetilde{\mathcal{N}}_{\widetilde{Q}} \in \mathcal{M}(\mathcal{U}^{\N})$:
$$
\widetilde{\mathcal{N}}_{\widetilde{Q}}(A)\, =\, \int_{\mathcal{M}_{\operatorname{Leb}}(\R\times \mathcal{U})} \widetilde{\nu}_{\gamma^{(1)}}(A)\, d\widetilde{Q}(\gamma)\, ,
$$
for each Borel set $A \subseteq \mathcal{U}^{\N}$.

\begin{theorem}
\label{thm:main}
(i) If $\widetilde{Q},\widetilde{Q}'$ are measures in $\mathcal{M}_{\text{st}}(\mathcal{M}_{\operatorname{Leb}}(\R\times \mathcal{U}))$ with $\widetilde{\mathcal{N}}_{\widetilde{Q}} = \widetilde{\mathcal{N}}_{\widetilde{Q}'}$,
then $\widetilde{Q}=\widetilde{Q}'$.\\[3pt]
(ii) The set of all thinning invariant measures $\mu \in \mathcal{M}(\mathcal{U}^{\N})$ 
is precisely 
$$
\{\widetilde{\mathcal{N}}_{\widetilde{Q}}\, :\, \widetilde{Q} \in \mathcal{M}_{\text{st}}(\mathcal{M}_{\operatorname{Leb}}(\R\times \mathcal{U}))\}\, .
$$
\end{theorem}
This is the most important theorem that we can prove, rigorously.
The proof is not difficult, but it does use a couple of nice ideas already in the literature.
In particular, we use a model introduced by Aldous,
and it seems that our proof is the first application of this process,
which Aldous introduced apparently for its sheer beauty.

\section{Proof of Theorem \ref{thm:main}}

Theorem \ref{thm:main} is proved by combining two ideas. The first is the idea of Poissonization
as in the paper of Ruzmaikina and Aizenman \cite{RuzmaikinaAizenman}.
The second is appeal to an important stochastic process that David Aldous described in his paper
\cite{Aldous} related to Hoyle's steady state model.
We describe the steady state model first.
This model is a Markov process version of the thinning invariant point processes, which were
proved to be just Poisson processes by Matthes, Kerstan and Mecke \cite{MatthesKerstanMecke}
and Kallenberg \cite{Kallenberg}.

Aldous's explanation of the steady state model recalls an old, now defunct cosmological model due
to Fred Hoyle, who sought a stationary alternative to the big bang theory.
Hoyle eventually converted to the big bang theory, himself.
See for example, \cite{Weinberg}.
But Aldous's description is a natural Markov process, equal to the time reversal of thinning.
We recall the thinning version first, which is easier.

We take the underlying space
to be $[0,\infty)$, and restrict attention to infinite, locally finite random point processes,
meaning that the random sample $\Xi \in \mathcal{M}_+([0,\infty))$, almost surely
may be expressed as 
$$
\Xi\, =\, \sum_{n=1}^{\infty}\delta_{R_n}\, ,
$$
for some sequence of points satisfying $0\leq R_1\leq R_2\leq \dots$, and $R_n \to \infty$.

Now, let $T_1,T_2,\dots$ in $[0,\infty)$ be i.i.d, random variables
$$
T_n\, \sim\, \operatorname{Exp}(1)\quad \text{ for each $n \in \N$,}
$$
and such that
$$
(T_1,T_2,\dots)\, \ind\, (R_1,R_2,\dots)\, .
$$
Then for each $t\geq 0$, let $\Xi_t \in \mathcal{M}_+([0,\infty))$
be the random measure
$$
\Xi_t\, =\, \sum_{n=1}^{\infty} \boldsymbol{1}_{[0,T_n)}(t)\, \delta_{R_n \exp(-t)}\, .
$$

If $\Xi$ is initially a homogeneous Poisson point process, e.g.,
$$
\Xi\, \sim\, \operatorname{PPP}([0,\infty),\|\cdot\|)\, ,
$$
then the process $\Xi_t$ is stationary. This is easy to see using the moment generating functional
from Definition \ref{def:MGF}.

For any $\tau \in \R$ we define the stochastic process $\Xi_{\tau,t}$ for $t \in [\tau,\infty)$
started at time $\tau$ as a Poisson point process on $[0,\infty)$ with intensity equal to Lebesgue measure
$\|\cdot\|$, and evolved according to the description above.
By Kolmogorov's extension theorem, these consistent distributions may be extended to 
give a law for a stationary process $\Xi_{-\infty,t}$ 
for all $t \in \R$.

Then one can consider the time reversed version, $(\Xi_{-\infty,-t})_{t \in \R}$, which is the version Aldous described, and which corresponds to Hoyle's steady state model.
At any given time $t$ there is a spatial Poisson point process of particles.
As one increases $t$ and considers $\Xi_{-\infty,-t}$ the particles move apart at a constant rate.
New particles are born into the universe to fill in the gaps.

\subsection{Marked version}
\label{subsubsec:marked}
Now, we consider the marked version of the steady state model.
Suppose that at time $\tau$ one has an initial random point process $\Xi$ now on $[0,\infty) \times \mathcal{U}$.
Let us suppose that taking the marginal on $[0,\infty)$, the point process $\Xi(\cdot \times \mathcal{U})$ is still
a Poisson process of intensity equal to Lebesgue measure.
Then almost surely there is a sequence $(R_1,U_1),(R_2,U_2),\dots \in [0,\infty) \times \mathcal{U}$,
such that
$$
\Xi\, =\, \sum_{n=1}^{\infty} \delta_{(R_n,U_n)}\, ,
$$
and
$$
0<R_1<R_2<\dots\, .
$$
Let $T_1,T_2,\dots$ be IID random variables
$$
T_1,T_2,\dots\, \sim\, \operatorname{Exp}(1)\quad \text{ and } \quad
(T_1,T_2,\dots)\, \ind\, ((R_1,U_1),(R_2,U_2),\dots)\, .
$$
Then for $t \in [\tau,\infty)$, define the random point process $\Xi_{\tau,t}$ such that for any Borel measurable
set $A \subseteq [0,\infty) \times \mathcal{U}$,
$$
\Xi_{\tau,t}(A)\, =\, \sum_{n=1}^{\infty} \boldsymbol{1}_{[0,T_n)}(t-\tau)\, \boldsymbol{1}_A(e^{\tau-t}R_n,U_n)\, .
$$
The marginal $\Xi_{\tau,t}(\cdot \times \mathcal{U})$ is the steady state model as before.

There are two important observations to remark upon.
Firstly, suppose that $\Xi=\Xi_{\tau,\tau}$ has some distribution as described above.
Let $U = (U_1,U_2,\dots)$.
Let us define $\mu$ in $\mathcal{M}(\mathcal{U}^{\N})$ to be the distribution
$$
\mu\, \stackrel{\text{def}}{:=}\, \mathscr{L}(U)\, .
$$
At any time $t \in [\tau,\infty)$, we may also almost surely find a sequence
$$
((R_1(t),U_1(t)),(R_2(t),U_2(t)),\dots)
$$ 
such that $0<R_1(t)<R_2(t)<\dots$ and
for any Borel set $A \subseteq [0,\infty) \times \mathcal{U}$,
$$
\Xi_{\tau,t}(A)\, =\, \sum_{n=1}^{\infty} \boldsymbol{1}_A(R_n(t),U_n(t))\, .
$$
Then, defining $U(t) = (U_1(t),U_2(t),\dots)$ and $\mu_{\tau,t}$ as 
$$
\mu_{\tau,t}\, \stackrel{\text{def}}{:=}\, \mathscr{L}(U(t))\, ,
$$
we have
$$
\mu_{\tau,t}\, =\, \Theta_{\exp(\tau-t)}(\mu)\, .
$$

Building on this, we make the second observation.
Let $\mathcal{M}_{\text{Leb}}([0,\infty) \times \mathcal{U})$ denote
the set of all Borel measures $\rho$ on $[0,\infty) \times \mathcal{U}$
such that
$\rho(\cdot \times \mathcal{U}) = |\cdot|$.
Given $\rho$ in this set, let us consider $\Xi$ to be the random point process
$$
\Xi\, \sim\, \operatorname{PPP}([0,\infty) \times \mathcal{U},\rho)\, .
$$
Then it satisfies the conditions above. More importantly, by the thinning property
of Poisson processes,
$$
\Xi_{\tau,t}\, \sim\, \operatorname{PPP}([0,\infty) \times \mathcal{U},\Phi_{t-\tau}(\rho))\, ,
$$
where for each $t \in [0,\infty)$, the mapping
$\Phi_{t} : \mathcal{M}([0,\infty)\times \mathcal{U}) \to \mathcal{M}([0,\infty)\times \mathcal{U})$ is non-random: let $\phi_t : [0,\infty) \to [0,\infty)$ be the mapping
$\phi_t(s) = st$, and then
$$
\Phi_t(\rho)\, =\, e^{-t} \rho \circ (\phi_{\exp(-t)} \times \operatorname{id}_{\mathcal{U}})^{-1}\, .
$$
By inspection, $\Phi_t$ sends the subset $\mathcal{M}_{\text{Leb}}([0,\infty)\times \mathcal{U})$ back to itself.

\subsection{Poissonization}

Given $\rho \in \mathcal{M}_{\text{Leb}}([0,\infty)\times \mathcal{U})$, let us denote a measure $\widehat{\nu}_{\rho}$ by taking
$$
\Xi\, \sim\, \operatorname{PPP}([0,\infty)\times \mathcal{U},\rho)\, ,
$$
which almost surely may be written as
$$
\Xi\, =\, \sum_{n=1}^{\infty} \delta_{(R_n,U_n)}\, ,
$$
for a sequence $(R_1,U_1),(R_2,U_2),\dots \in [0,\infty) \times \mathcal{U}$ with $0<R_1<R_2<\dots$, and then defining
$$
\widehat{\nu}_{\rho}\, \stackrel{\text{def}}{:=}\, \mathscr{L}(U_1,U_2,\dots)\, .
$$
Given $Q \in \mathcal{M}\left(\mathcal{M}_{\operatorname{Leb}}([0,\infty)\times \mathcal{U})\right)$,
we define $\widehat{\mathcal{N}}_Q$ such that
$$
\widehat{\mathcal{N}}_Q(A)\, =\, \int_{\mathcal{M}_{\operatorname{Leb}}([0,\infty)\times \mathcal{U})} \widehat{\nu}_{\rho}(A)\, dQ(\rho)\, ,
$$
for every Borel set $A\subseteq \mathcal{U}^{\N}$.

\begin{lemma}
\label{lem:divis}
Suppose that $\mu \in \mathcal{M}(\mathcal{U}^{\N})$ is thinning invariant. Then
there is a Borel probability measure  
$Q \in \mathcal{M}\left(\mathcal{M}_{\operatorname{Leb}}([0,\infty)\times \mathcal{U})\right)$,
such that $\mu = \widehat{\mathcal{N}}_Q$.
\end{lemma}
The proof is simple. We use an idea from \cite{RuzmaikinaAizenman}.
The idea of Poissonization is natural in these problems.

\begin{proof}
Suppose $\mu$ is thinning invariant.
Let $U=(U_1,U_2,\dots)$ be random and $U \sim \mu$.
Independently, let $(R_1,R_2,\dots)$ be distributed according to a standard
Poisson point process on $[0,\infty)$. Take
$$
\Xi\, \stackrel{\text{def}}{:=}\, \sum_{n=1}^{\infty} \delta_{(R_n,U_n)}\, .
$$
This is not Poissonian, but we will push it back in time to $-\infty$ to obtain something asymptotically
Poissonian at time $0$.

Consider the following approximation.
Let $\Xi'$ be a random point process, such that, conditional on $\Xi$, we have
$$
\Xi'\, \sim\, \operatorname{PPP}([0,\infty)\times \mathcal{U},\Xi)\, .
$$
Then for any $\tau \in \R$ and $t \in [\tau,\infty)$, the process described before
as $\Xi'_{\tau,t}$ is essentially the same as $\Xi_{\tau,t}$ with the following difference:
In $\Xi_{\tau,t}$ each original point receives a Bernoulli random variable with survival rate
$p = \exp(\tau-t)$ to determine if it remains or is deleted;
whereas, in $\Xi'_{\tau,t}$, each original point is copied a Poisson number of times,
with rate $p=\exp(\tau-t)$.
We want to compare $\Xi_{\tau,0}$ and $\Xi'_{\tau,0}$ in the limit where $\tau \to -\infty$.

Suppose $n$ is a fixed integer.
If we keep track of the first $n/p$ of the original points of $\Xi$, whether we act on them according 
to Bernoulli thinning or ``Poisson thinning'' we see on average $n$ survivors which is an order-1
quantity. More precisely in the first setting the number of survivors is binomial with parameters
$(n/p,p)$ and in the latter case it is Poisson with parameter $n$.
By the usual law of small numbers the number of survivors converge in distribution in the limit $p \to 0$.
Note that if we take $t=0$ then $p = \exp(\tau)$ so that $p \to 0$ as $\tau \to -\infty$.

In addition to just keeping track  of the number of points, we should also check the point-values
themselves.
According to the weak topology, if we can couple each of these processes with a small
failure rate, then that suffices to show convergence to zero of the difference in the limit as $\tau \to -\infty$.
The failure rate for the Poisson-to-Bernoulli coupling is proportional to the probability of having two or more survivors
according to any of the $n/p$ original Poisson variables.
Since each Poisson random variable has rate $p$, and since the probability of two
or more points is on the order of $p^2$, this shows that the overall failure rate is bounded by a constant times $np$. This converges to zero
for each fixed $n \in \N$, in the limit as $p \to 0$.

Since we have shown that the distributions of $\Xi_{\tau,0}$ and $\Xi'_{\tau,0}$ are asymptotically
close in the limit $\tau \to -\infty$, we now want to check that there is a limit-point for the distribution $\Xi'_{\tau,0}$ which has the desired properties.
One can take the random meaure $\rho=\Xi$ so that, conditional on this value
$$
\Xi'\, \sim\, \operatorname{PPP}([0,\infty)\times \mathcal{U},\rho)\, .
$$ 
Note that, conditioning on the random value of $\rho=\Xi$, we have
$$
\Xi'_{\tau,0}\, \sim\, \operatorname{PPP}([0,\infty)\times \mathcal{U},\Phi_{-\tau}(\rho))\, .
$$
Now by the strong law of large numbers, one can see that for almost all choices of $\rho$,
$$
\lim_{\tau \to -\infty} \Phi_{-\tau}(\rho)(A \times \mathcal{U})\, =\, |A|\, .
$$
This is just a fact about the original steady state model because we are taking the marginal
on $[0,\infty)$. For instance $\Phi_{-\tau}(\rho)([0,x)\times \mathcal{U})$ just measures
$$
e^{\tau} |\{n\, :\, R_n < x \exp(-\tau)\}|\, ,
$$ 
for the original Poisson point process $(R_1,R_2,\dots)$
on $[0,\infty)$.
From this, we see that if there is a limit-point for the distributions of $\Phi_{-\tau}(\rho)$,
along some subsequence $(\tau_1,\tau_2,\dots)$ with $\tau_n \to -\infty$,
then this limit will be a Borel probability measure $Q$ which is supported on 
$\mathcal{M}_{\text{Leb}}([0,\infty)\times \mathcal{U})$.

To guarantee that a subsequential limit point exists, we just need to check tightness. 
(See, for example, \cite{Durrett}.)
One may break up $[0,\infty)$ into countably many intervals $[n-1,n)$
for $n \in \N$.
Thus we may consider $\mathcal{M}_+([0,\infty)\times \mathcal{U})$ as
$$
\prod_{n \in \N} \mathcal{M}_+([n-1,n)\times \mathcal{U})\, ,
$$
and in fact the product topology is appropriate for the local weak topology on point processes
on $[0,\infty) \times \mathcal{U}$.
In the last paragraph we already established convergence to a Poisson point process for the marginal on $[0,\infty)$, which guarantees tightness
for the marginal on each interval $[n-1,n)$. Since $\mathcal{U}$ is compact the marginal
on $\mathcal{U}$ is automatically tight. Also the product topology of any number of compact sets
is also compact. This guarantees tightness overall.
\end{proof}

\subsection{Proof of Theorem \ref{thm:main} part (i)}

Suppose $\mu \in \mathcal{M}(\mathcal{U}^{\N})$ is thinning invariant. Then, according to Lemma \ref{lem:divis}, there is some Borel probability measure $Q$ supported on $\mathcal{M}_{\text{Leb}}([0,\infty)\times \mathcal{U})$
such that $\mu = \widehat{\mathcal{N}}_Q$.
By the comments at the end of Subsection \ref{subsubsec:marked}, it is clear that
for any $\tau<0$, we could replace $Q$ by $\overline{Q}_{[\tau,0]}$, the ergodic average
$$
\overline{Q}_{[\tau,0]}(\cdot)\, \stackrel{\text{def}}{:=}\, \frac{1}{-\tau} \int_{\tau}^{0} Q\circ \Phi_t^{-1}(\cdot)\, dt\, .
$$
One still has the desired condition
that $\mu = \widehat{\mathcal{N}}_{\overline{Q}_{[\tau,0]}}$

The total variation norm of $(\overline{Q}_{[\tau,0]} \circ \Phi_t^{-1} - \overline{Q}_{[\tau,0]})$ is bounded by $t/|\tau|$.
By tightness/compactness, we may find a sequence $(\tau_1,\tau_2,\dots)$ with $\tau_n \to -\infty$
such that $\overline{Q}_{[\tau_n,0]}$ converges.
Therefore, we obtain $Q^*$, a Borel probability measure supported on $\mathcal{M}_{\text{Leb}}([0,\infty)\times \mathcal{U})$ such that $Q^* \circ \Phi_t^{-1} = Q^*$ for every $t \in [0,\infty)$,
and such that $\mu = \mathcal{N}_{Q^*}$.
For notational convenience, we will assume that this measure was chosen at the outset: $Q=Q^*$.

With this, we can set-up a steady state model on $[0,\infty) \times \mathcal{U}$ with
stationary measure given by first selecting a density $\rho \in \mathcal{M}_{\text{Leb}}([0,\infty)\times \mathcal{U})$ distributed according to $Q$, and then taking the random point process $\Xi$,
such that conditional on $\rho$,
$$
\Xi\, \sim\, \operatorname{PPP}([0,\infty)\times \mathcal{U},\rho)\, .
$$
Note that for $\tau \in \R$ and $t \in [\tau,\infty)$, conditional on $\rho$, we have
$$
\Xi_{\tau,t}\, \sim\, \operatorname{PPP}([0,\infty)\times \mathcal{U},\Phi_{t-\tau}(\rho))\, .
$$
But since $\rho \sim Q$ and since $Q \circ \Phi_t^{-1} = Q$ for all $t \in [0,\infty)$, this means
we still have
$$
\Phi_{t-\tau}(\rho)\, \sim\, Q\, .
$$
So the overall distribution, averaging over $\rho$, is invariant.
By Kolmogorov's extension principle, we can now extend to obtain the model when $\tau \to -\infty$.

\subsubsection{Uniqueness}
Suppose that $\Xi_{-\infty,t}$ is this process. In other words, at each time $t$, it is a point process
on $[0,\infty) \times \mathcal{U}$: its marginal distribution at any given time is as a Cox process
with mixing measure $Q$, and it evolves as $t$ increases as indicated above.
We may define $\mathcal{F}_{(-\infty,\tau]}$ to be the $\sigma$-algebra of $(\Xi_{-\infty,t})_{t\leq \tau}$.
Then, for any measurable set $A \subseteq [0,\infty) \times \mathcal{U}$, we may
consider the reversed (or backwards) Doob's martingale
$$
\E[\Xi_{-\infty,0}(A)\, |\, \mathcal{F}_{(-\infty,-t]}]\quad \text{ for } \quad t \in [0,\infty)\, .
$$
According to Doob's reversed martingale convergence theorem, we know that this converges, a.s,
and in expectation.
(See for example \cite{Durrett}, Section 5.6.)
Let us call the limit $\widetilde{\rho}(A)$. Note that this is measurable with respect to $\mathcal{F}_{-\infty}$, the backwards tail $\sigma$-algebra.
On the other hand, using the notation $(R_n(t),U_n(t))_{n=1}^{\infty}$ introduced before, we see that
for $t \in [0,\infty)$ the conditional distribution of $\Xi_{-\infty,0}(A)$ given $\mathcal{F}_{(-\infty,-t]}$ is 
\begin{equation}
\label{eq:twopass}
\E\left[e^{\displaystyle -\lambda \Xi_{-\infty,0}(A)}\, \Big|\, \mathcal{F}_{(-\infty,-t]}\right]\,
=\, e^{\displaystyle e^t \Phi_t( \Xi_{-\infty,-t})(A) \ln \left(1-(1-e^{-\lambda})e^{-t}\right)}\, .
\end{equation}
Using (\ref{eq:twopass}) in a first pass, 
one may see directly that
$$
\E[\Xi_{-\infty,0}(A)\, |\, \mathcal{F}_{(-\infty,-t]}]\, =\, \Phi_t( \Xi_{-\infty,-t})(A)\, .
$$
Using this fact, and making a second pass at (\ref{eq:twopass}), we see that,
conditional on $\mathcal{F}_{-\infty}$,
$\Xi_{-\infty,0}(A)$ is distributed as a Poisson random variable with parameter $\widetilde{\rho}(A)$.

To gain insight into the distribution of $\widetilde{\rho}$, note that by stationarity of $\Xi_{-\infty,t}$,
the random measure $\Phi_t(\Xi_{-\infty,-t})$
has the same distribution as $\Phi_t(\Xi_{-\infty,0})$.
But the distribution of $\Xi_{-\infty,0}$ is obtained as follows: let $\rho$ be random with $\rho \sim Q$,
and conditional on this $\Xi_{-\infty,0} \sim \operatorname{PPP}([0,\infty) \times \mathcal{U},\rho)$.
Moreover, since $Q \circ \Phi_t^{-1} = Q$, we may make yet another change and still have the same distribution. Choose $\rho$ randomly with $\rho \sim Q$, and then conditionally on this, for each $t \in [0,\infty)$, let $\widetilde{\Xi}_t$ be random with
$$
\widetilde{\Xi}_t\, \sim\, \operatorname{PPP}([0,\infty)\times \mathcal{U},\Phi_{-t}(\rho))\, .
$$
Note that $\Phi_t$ is well defined for $t<0$ as well as for $t\geq 0$.

Then the overall distribution of $\Phi_t(\widetilde{\Xi}_t)$ is still the same as $\Phi_t(\Xi_{-\infty,t})$,
after averaging over $\rho$ as well as the conditional Poisson measure, given $\rho$.
On the other hand, we claim that, conditional on $\rho$, we almost surely have convergence $\Phi_t(\widetilde{\Xi}_t) \to \rho$ in distribution.
To see this,  note that for any $f \in \mathcal{B}_+([0,\infty)\times \mathcal{U})$, 
$$
\E\Bigg[e^{\displaystyle -\int_{[0,\infty) \times \mathcal{U}} f d\Phi_t(\Xi_{-\infty,t})}\Bigg]\,
=\, \exp\left(\int_{[0,\infty)\times \mathcal{U}} \left(1-\exp\left(-e^{-t} f\right)\right) e^t\, d\rho\right)\, ,
$$
using Definition \ref{def:MGF} and the definition of $\Phi_t$. By the monotone convergence theorem,
this converges to $\exp(-\int f\, d\rho)$ which is the moment generating function for $\rho$,
and this uniquely characterizes the random measure.
This implies that the distribution of $\widetilde{\rho}$ is $Q$.
Since $\widetilde{\rho}$ was obtained from the process, itself, using the martingale convergence theorem,
this is the usual proof of uniqueness of $Q$ (similar to the proof of uniqueness for the directing measure of infinite
exchangeable sequences of random variables in de Finetti's theorem \cite{AldousReview}).

\subsection{Births and a change of variables}
\label{sec:birth}

We now consider the steady state model we have been considering in reversed time $(\Xi_{-\infty,-t})_{t \in \R}$.
This was the perspective taken by Aldous.
In this perspective new points are being created by a space-time Poisson process with marks.
In the thinning perspective these would be extinction times for particles coming from $-\infty$.
From this perspective, 
$\Xi_{-\infty,-t}(A) = H(\Psi_t(A))$, where $H$ is a space-time-marking point process and
$$
\Psi_t(A)\, =\, \{(s,re^{-t-s},u)\, :\, (r,u) \in A\, ,\ s \in [-t,\infty)\}\, .
$$
One may determine $H$ by starting with the distribution of $\Xi_{-\infty,t}$
and then performing deletions, and keeping track of the points as they are deleted.

At a time $t$, the distribution of points is $\Xi_{-\infty,t}$ which is Poisson on $[0,\infty)\times \mathcal{U}$,
with intensity $\Phi_{t}(\widetilde{\rho})$, where we are conditioning on $\mathcal{F}_{-\infty}$,
continuing with the analysis of the last section.
Each such point is deleted at uniform rate. So in an infinitesimal time interval $[t,t+dt)$,
the distribution of deleted points is $\Phi_{t}(\widetilde{\rho})(dr\times du)\, dt$.
But a Poisson process has the thinning property: separating a Poisson point process
into two groups according to an independent Bernoulli process both groups are independent
Poisson processes with certain rates.
Using this and the infinitesimal generator of thinning
$\Phi_{t}(\widetilde{\rho})(dr\times du)\, dt$,
we see
$$
H\sim \operatorname{PPP}(\R \times [0,\infty) \times \mathcal{U},\widehat{\Phi}(\widetilde{\rho}))\quad \text{ where} \quad
\widehat{\Phi}(\widetilde{\rho})(dt\times dr\times du)\, =\, \Phi_{t}(\widetilde{\rho})(dr\times du)\, dt\, .
$$
Let us rewrite this, using the definition of $\Phi_t$:
$$
\widehat{\Phi}(\widetilde{\rho})([t_1,t_2)\times[r_1,r_2)\times A)\,
=\, \int_{t_1}^{t_2} e^{-t} \widetilde{\rho}([e^{t}r_1,e^{t}r_2)\times A)\, dt\, .
$$
Now we will do a change of variables. The calculations are direct, but not transparent,
so we will perform a consistency check. 
We will recover the marginal
density $\Xi_{-\infty,t}$ from $H$ using the new variables.

Consider the function $\mathcal{L}(t,r) = (t,\ln(r)+t)$ from $\R \times (0,\infty)$ to $\R \times \R$.
We may consider rectangles
$$
\widehat{\Phi}(\widetilde{\rho})(\mathcal{L}^{-1}([t_1,t_2)\times [x_1,x_2)) \times B)\,
=\, \int_{t_1}^{t_2} e^{-t} \widetilde{\rho}([\exp(x_1),\exp(x_2)) \times B)\, dt\, ,
$$
for any Borel set $B \subseteq \mathcal{U}$.
Note that the inverse mapping is $\mathcal{E}(t,x) = (t,e^{x-t})$.
Therefore,
$$
\Psi_t([r_1,r_2)\times B)\, =\, \mathcal{L}^{-1}([-t,\infty) \times [\ln(r_1)-t,\ln(r_2)-t)) \times B\, .
$$
Hence, the expectation of
$\Xi_{-\infty,-t}([r_1,r_2)\times B) = H(\Psi_t([r_1,r_2) \times B))$ is
$$
\widehat{\Phi}(\widetilde{\rho})(\Psi_t([r_1,r_2)\times B))\,
=\, \int_{-t}^{\infty} e^{-s}\, \widetilde{\rho}([r_1e^t,r_2e^t)\times B)\, ds\, 
=\, e^t \widetilde{\rho}([r_1e^t,r_2e^t)\times B)\, ,
$$
which is $\Phi_{-t}(\widetilde{\rho})([r_1,r_2)\times B)$.
This is as it should
be since $\Phi_{-t}(\widetilde{\rho})$ is the intensity (conditional on $\mathcal{F}_{-\infty}$)
for the Poisson point process $\Xi_{-\infty,-t}$,
given $\widetilde{\rho}$. This is a consistency check that we have applied the inverse functions $\mathcal{L}(t,r)$
and $\mathcal{E}(t,x)$ correctly to the sets and the measures.

\subsubsection{Proof of Theorem \ref{thm:main} part (ii)}

Now let us consider $H \circ (\mathcal{L} \times \operatorname{id}_{\mathcal{U}})^{-1}$.
If we define $\gamma^{(1)} \in \mathcal{M}_+(\R \times \mathcal{U})$ by
$$
\gamma^{(1)}([x_1,x_2)\times B)\, =\, \widetilde{\rho}([\exp(x_1),\exp(x_2))\times B)\, ,
$$
then notice that $\gamma^{(1)}([x_1,x_2)\times \mathcal{U}) = \int_{x_1}^{x_2} e^x\, dx$.
Therefore, 
$$
\gamma^{(1)}(dx\times du)\, =\, e^x\, \gamma(dx\times du)\, ,
$$ 
where $\gamma$
is in $\mathcal{M}_{\text{Leb}}(\R\times \mathcal{U})$.
The calculations of the last subsection show that
$$
H \circ (\mathcal{L}\times \operatorname{id}_{\mathcal{U}})^{-1}\,
\sim\, \operatorname{PPP}(\R\times \R\times \mathcal{U},\eta \times \gamma)\, ,
$$
where $\eta$ is the meaure on $\R$: $d\eta(t) = e^{-t}\, dt$.

Using the fact that $\Xi_{-\infty,-t}([r_1,r_2)\times B) = H(\Psi_t([r_1,r_2)\times B))$,
we see that stationarity of $(\Xi_{-\infty,-t})_{t \in \R}$ amounts to stationarity in $t$ of the marginal distribution of 
\begin{equation}
\label{eq:stat2}
H\circ \mathcal{L}^{-1}([-t,\infty)\times [x_1-t,x_2-t) \times B)\, ,
\end{equation}
averaging over both the conditional Poisson distribution given $\gamma$, and over the induced
measure on $\gamma$.
Let us define $\widetilde{Q}$ to be the distribution of the random measure $\gamma$ coming from $Q$ 
by pulling back the transformation $\widetilde{\rho} \mapsto \gamma$.
Then thinning invariance at the level $H' = H \circ (\mathcal{L} \times \operatorname{id}_{\mathcal{U}})^{-1}$ as in (\ref{eq:stat2})
is stationarity with respect to the
shift: $\widetilde{\sigma}_t : \R \times \R$ given by $\widetilde{\sigma}_t(s,x) = (s+t,x+t)$. Note that the
Poisson point process $H'$ has an intensity which may written using $\gamma$:
$$
d\eta(s) \times \gamma^{(1)}(dx\times du)\, =\, e^{x-s}\, \gamma(dx\times du)\, ds\, .
$$
This satisfies the property that 
$$
(\eta \times \gamma^{(1)}) \circ (\widetilde{\sigma}_t \circ \operatorname{id}_{\mathcal{U}})^{-1}\,
=\, \eta \times [\gamma \circ (\sigma_t \times \operatorname{id}_{\mathcal{U}})^{-1}]^{(1)}\, .
$$
Therefore, stationarity of $(\Xi_{-\infty,-t})_{t \in \R}$ actually amounts to stationarity of the Cox process on $\R \times \mathcal{U}$,
relative to the usual shifts $(\sigma_x \times \operatorname{id}_{\mathcal{U}})$, where the
the directing measure of the Cox process is $\widetilde{Q}$.
But this is equivalent to stationarity of the directing measure $\widetilde{Q}$.
See for example, Daley and Vere-Jones \cite{DaleyVere-Jones}, Section 10.1, Exercise 10.1.3,
which is an easy exercise in using the definition of the moment generating functional as in Definition \ref{def:MGF}.
This is the characterization that was to be proved.
So part (ii) of the theorem is completed.

\section{Thinning for Gap Distributions}
\label{sec:gap}

Let $\Gamma \subset [0,\infty)^{\N}$ consist of all sequences $w=(w_1,w_2,\dots)$ satisfying
$$
0\, \leq\,  w_1\, \leq\,  w_2\, \leq\,  \dots\, .
$$
We call these configurations of gaps, for a reason which will be apparent momentarily.

Given $\mu \in \mathcal{M}(\Gamma)$ and $p \in (0,1]$ we consider a new measure
$\widehat{\theta}_p(\mu)$ in $\mathcal{M}(\Gamma)$ defined as follows.
Let $W = (W_1,W_2,\dots)$ be a random element of $\Gamma$,
$$
W\, \sim\, \mu\, .
$$
Define $W_0 = 0$.
Additionally let $B = (B_0,B_1,\dots)$ be i.i.d., Bernoulli-$p$ random variables
$$
B\, \ind\, W\, .
$$
Let $K_0<K_2<\dots$ be the subsequence of $\{0,1,\dots\}$ uniquely defined, a.s., by the condition
$$
\{K_0,K_1,\dots\}\, =\, \{k \in \{0,1,\dots\}\, :\, B_k=1\}\, .
$$
Let $V_0=0$ and $V_n=(W_{K_n}-W_{K_0})$ for $n \in \N$.
These are gaps of the thinned sequence starting from $(W_0,W_1,\dots)$. Let us define $V = (V_1,V_2,\dots)$,
which is an element of $\Gamma$, a.s.
Then
$$
\widehat{\theta}_p(\mu)\, \stackrel{\text{def}}{:=}\, \mathscr{L}(V)\, .
$$

\subsection{Poisson construction for gap distributions}

Let $\mathcal{M}_{\uparrow}(\R)$ denote the set of all measures $\alpha \in \mathcal{M}_{+}(\R)$, satisfying
\begin{itemize}
\item $\alpha((-\infty,x))<\infty$ for each $x \in \R$, and 
\item $\alpha(\R) = \infty$.
\end{itemize}
Given $\alpha \in \mathcal{M}_{\uparrow}(\R)$ we define $\widehat{\nu}_{\alpha} \in \mathcal{M}(\Gamma)$ as follows.
Take $\Xi$ a random point process on $\R$ with
$$
\Xi\, \sim\, \operatorname{PPP}(\R,\alpha)\, .
$$
Almost surely, there is a sequence $-\infty<X_1\leq X_2\leq \dots$, with
$$
\Xi\, =\, \sum_{n=1}^{\infty} \delta_{X_n}\, .
$$
Let $U_n = (X_{n+1}-X_n)$ for $n \in \N$, and note that $U = (U_1,U_2,\dots)$ is an element of 
$\Gamma$, a.s.
We define
$$
\widehat{\nu}_{\alpha}\, \stackrel{\text{def}}{:=}\, \mathscr{L}(U)\, .
$$
It is easy to see that $\widehat{\nu}_{\alpha \circ \sigma_x^{-1}} = \widehat{\nu}_{\alpha}$
fo all $x \in \R$.
Also,
$$
\widehat{\theta}_p(\widehat{\nu}_{\alpha})\, =\, \widehat{\nu}_{p\alpha}\, .
$$
Note that (modulo shifts) the only measures $\alpha \in \mathcal{M}_{\uparrow}(\R)$ satisfying
$$
p \alpha\, =\, \alpha \circ \sigma_{x(p)}^{-1}\, ,
$$
are $\alpha_m(dx) = m e^{mx}\, dx$ for $m \in (0,\infty)$, which each leads
to a thinning invariant gap distribution $\widehat{\nu}_{\alpha_m}$.
But the limit 
$\lim_{m \downarrow 0} \widehat{\nu}_{\alpha_m}$ also exists, and is the point mass
at $(0,0,\dots) \in \Gamma$. We denote this as $\widehat{\nu}_{\infty \delta_0}$.
These are the ``pure'' Poisson gap distributions which are thinning invariant.

\subsection{Cox construction and conjecture for gaps}
While we have classified the ``pure'' Poisson gap distributions which are thinning invariant,
one may generically perturb each of these by some version of stationary Cox processes
to obtain something which is still thinning invariant, in complete analogy to the case of
thinning invariant sequences.

Let $\mathscr{K}'$ denote the disjoint union of $\mathcal{M}_{+,1}(\R) \times (0,\infty)$ and $\{0\}$.
Given $(\gamma,m) \in \mathcal{M}_{+,1}(\R) \times (0,\infty)$, we define
$$
\gamma^{(m)}(dx)\, \stackrel{\text{def}}{:=}\, m e^{mx}\, \gamma(dx)\, .
$$
(Note that this is consistent with the definition (\ref{eq:1prime}) in the case $m=1$.)
We define $\mathcal{M}_{\text{st}}(\mathscr{K}')$ to be the set of all measures $Q \in \mathcal{M}(\mathscr{K}')$
such that for each Borel set $A \subseteq \mathcal{M}_{+,1}(\R) \times (0,\infty)$ and each $x \in \R$
$$
Q(\{(\gamma,m)\, :\, (\Sigma_x(\gamma),m) \in A\})\, =\, Q(A)\, .
$$
Given $Q \in \mathcal{M}_{\text{st}}(\mathscr{K})$ we define $\widehat{\mathcal{N}}_Q \in \mathcal{M}(\Gamma)$
such that for any Borel set $A \subseteq \Gamma$,
$$
\widehat{\mathcal{N}}_{Q}(A)\, =\, Q(\{0\}) \widehat{\nu}_{\infty \delta_0}(A) + 
\int_{\mathcal{M}_{+,1}(\R) \times (0,\infty)} \widehat{\nu}_{\gamma^{(m)}}(A)\, dQ(\gamma,m)\, .
$$

\begin{conjecture}
\label{conj:gap}
(i) For any pair $Q,Q' \in \mathcal{M}_{\text{st}}(\Gamma)$ if $\widehat{\mathcal{N}}_Q = \widehat{\mathcal{N}}_{Q'}$
then $Q=Q'$.\\[5pt]
(ii) The set of thinning invariant $\mu \in \mathcal{M}(\Gamma)$ is $\{\widehat{\mathcal{N}}_Q\, :\, Q \in \mathcal{M}_{\text{st}}(\Gamma)\}$.
\end{conjecture}

\subsection{Heuristic argument for gaps}

We will outline an argument for Conjecture \ref{conj:gap}, paying special attention to the most important step that we cannot rigorously prove.

Suppose that $\mu$ is a thinning-invariant gap distribution. We may choose $W \sim \mu$
where $W = (W_1,W_2,\dots)$.
We define a point process $\Xi_0 = \sum_{n=1}^{\infty} \delta_{W_n}$.
Now, let $T_1,T_2,\dots$ be independent $\operatorname{Exp}(1)$ random variables
$$
(T_1,T_2,\dots) \ind (W_1,W_2,\dots)\, ,
$$
and define $\Xi_t = \sum_{n=1}^{\infty} \boldsymbol{1}_{[0,T_n)}(t)\, \delta_{W_n}$.
Then modulo a shift this is the thinning of $\Xi_0$.
Therefore, there is a random constant $X(t)$ such that $\Xi_t\circ \sigma_{X_t}^{-1}$
is stationary.
From the semi-group property of deletions, it is easy to see that $X_t$ itself is stationary.
Therefore, by the ergodic theorem, there is a random constant $M$ such that
\begin{equation}
\label{eq:ergodic}
\lim_{t \to \infty} \frac{X_t}{t}\, =\, 1/M\, .
\end{equation}
(See for example, Chapter 7 of  \cite{Durrett}.)
In order to simplify the situation, we will suppose that $\mu$ satisfies the property that $M$ is non-random,
$\mu$-almost-surely, and consider this value to be $m$.
We then shift the originally defined point processes:
$$
\Xi'_t\, =\, \Xi_t \circ \sigma_{t/m}^{-1}\, .
$$
We have taken care of the leading order effect of the shift.
This is where we need to make a non-rigorous jump.

\begin{conjecture}
(i) The case $M=0$ is not possible in (\ref{eq:ergodic}).\\[3pt]
(ii) For $M\in(0,\infty]$, define $(N_1(t),N_2(t),\dots)$ such that $N_1(t)<N_2(t)<\dots$ and $(N_k(t))_{k=1}^{\infty}  = (n\, :\, T_n<t)$ for each $t$.
Then the distribution of $(W_{N_k(t)} - m^{-1} t)$ is tight in $t$ for each $k \in \N$.
\end{conjecture}

Let us consider $m \in (0,\infty)$ first, before considering $m=\infty$.
We just quickly repeat the steps from the proof of Theorem \ref{thm:main}, making the appropriate
changes for the present situation.

We may start at time $\tau$: $\Xi'_{\tau,t}$ is $\Xi_{t-\tau}$ for $t>\tau$.
Then by tightness, 
we may take a subsequence $\tau_n$ with $\tau_n \to -\infty$ such that the weak limit of the entire
distribution of $\Xi'_{\tau_n,t}$ converges for all $t \in \R$.
Then we use Kolmogorov's extension theorem to construct $\Xi'_{-\infty,t}$.
Then we may take the ergodic/Cesaro limit to obtain a stationary version which we assume we had done from the outset
to avoid complicating the notation.
(Note that shifting by any times does not affect the marginal distribution of gaps under the thinning-invariant assumption.)
One can then consider the process in reverse.

Now, conditional on the backwards tail-algebra $\mathcal{F}_{-\infty}$,
particles are born (backwards in time) according to a space-time process $H \sim \operatorname{PPP}(\R \times \R, \hat{\rho}\circ \sigma_{t/m}^{-1}(dx) dt)$
where $\hat{\rho}$ is obtained by keeping track of particles in space-time, streaming at the constant velocity $1/m$.
After birth the particles drift to the right at this velocity $1/m$.
Stationarity means that the new particles are filling in gaps caused by the drift.
Although the drift is constant in space, therefore preserving Lebesgue measure, recall that the pure Poisson construction
for parameter $m$ is based on a Poisson process with intensity $m e^{mx}$. So that the particles are more dense to the right,
which means if the entire point process drifts to the right, then the gaps become relatively bigger (since they are the gaps that were previously on the left).
We claim that stationarity of $\Xi'_t$ implies that $\hat{\rho} = \gamma^{(m)}$ for $\gamma$ stationary (with respect to shifts).
We leave this as an exercise similar to Section \ref{sec:birth}, rewriting stationarity of $\Xi'_t$ in terms of properties
of the Cox process $H$, and then using the exercise from Daley and Vere-Jones to relate stationarity of the Cox process
to stationarity of its directing measure.

Notice that for $m \in (0,\infty)$, equation (\ref{eq:ergodic}) should imply that $W_n-W_0/\ln n \to 1/m$.
One can see this by noting that $N_1(t)$ is geometrically distributed with failure rate $p=e^{-t}$.
Therefore this claim follows by a suitable Tauberian theorem to de-Geometrize.
We would like to refer the reader to \cite{Berestycki}, Section 1.5, for important uses of Tauberian theorems
in the context of coalescents, which are illuminating.

If $m=0$ then tightness of $(W_{N_1(t)})_{t \in \R}$ implies tightness of $(W_1,W_2,\dots)$
(using another Tauberian argument, as is necessary).
But due to monotonicity this implies there is a limsup, and then thinning-invariance implies that $W=(0,0,\dots)$, a.s.

We now make an important remark.

\begin{remark}
The conjecture should not be trivial. Our proof of Theorem \ref{thm:main} used
one idea from the paper of Ruzmaikina and Aizenman \cite{RuzmaikinaAizenman}: Poissonization.
But their characterization of gap distributions invariant under uncorrelated cavity steps involves
other arguments whose analogues we have not yet found.
In principle the present problem should require an equal amount of work.
\end{remark}

\section{Return to partition structures}

\label{sec:return}

One can move from partition structures to gap distributions by replacing $\xi = (\xi_1 \geq \xi_2 \geq \dots)$
with $W = (W_1,W_2,\dots)$ with $W_n = \ln\xi_1-\ln\xi_{n+1}$ for $n \in \N$.
In thinning the partition structure one rescales $\xi$.
At the level of $W$ this amounts to a particular random shift of $(-\ln\xi_1,-\ln\xi_2,\dots)$.
But the gaps are defined from the perspective of the leading point (which is $\ln\xi_1$, but we find it convenient
to introduce a reflection).
So any thinning invariant partition structure leads to a thinning invariant gap distribution.

Note that the proportion of dust in a random partition structure is non-decreasing under $\theta_p$
and it is constant only if the total amount of dust is $0$ or $1$.
Therefore a thinning invariant random partition structure has an amount of dust $\{0,1\}$, a.s.
Full dust is $\nu_{(0,1)}$.
Otherwise, with no dust, the sequence $W$ can be used to fully recover $\xi$, where $\xi_1$ is determined by the normalization condition.

We can then use the conjectured characterization of the thinning invariant gap distributions
to find the thinning invariant random partition structures.
But due to the normalization condition this requires $m \in (0,1)$.
The limit as $m \to 1$ leads to $\nu_{(0,1)}$.
With this, Conjecture \ref{conj:main} follows from Conjecture \ref{conj:gap}.

\subsection{Relation to REM}

Let $B_n(t)$ be independent Brownian motions for $n \in \N$, centered (and with variance equal to $t$ at time $t$).
Consider weights $X_{N,n}(\beta) = e^{-\beta B_n(N)}$ for $n=1,\dots,2^N$.
Then the REM partition function is 
$$
Z_N(\beta)\, =\, \frac{1}{N}\, \ln \sum_{n=1}^{2^N} X_{N,n}(\beta)\, .
$$
One observation is that for $n \in \{1,\dots,2^N\}$, we have
$$
\ln X_{N+1,n}(\beta)\, =\, -\beta [B_n(N+1)-B_n(N)] + \ln X_{N,n}(\beta)\, ,
$$
and $[B_n(N+1)-B_n(N)]$ is independent of $X_{N,n}$.
This makes contact with the uncorrelated cavity step dynamics introduced by Aizenman and Ruzmaikina.
But additionally, we have $2^N$ new points at time step $N+1$: $X_{N+1,n}(\beta)$ for $n=2^N+1,\dots,2^{N+1}$,
and we mentioned in the introduction how this leads to the notion of thinning.

Suppose that $\xi \in \Delta_{\infty}$ is a sample of a  random partition structure.
Let us define two random constants:
$$
\mathscr{E}(\xi)\, \stackrel{\text{def}}{:=}\, -\frac{1}{\beta}\, \ln\left( \sum_{n=1}^{\infty} \xi_n e^{-\beta Z_n} + \xi_0 \E[e^{-\beta Z_0}]\right)\, ,
$$
where $Z_0,Z_1,\dots$ are independent $\mathcal{N}(0,1)$ random variables
$$
(Z_0,Z_1,\dots)\, \ind\, \xi\, ,
$$
and
$$
\mathscr{S}(\xi)\, =\, \ln\left( \sum_{n=1}^{\infty} \xi_n B_n \right)\, ,
$$
where $B_1,B_2,\dots$ are independent $\operatorname{Bernoulli}(1/2)$ random variables independent of $\xi$.
Then one would expect that the free energy of the REM
$$
\mathscr{F}(\beta)\, \stackrel{\text{def}}{:=}\, \lim_{N \to \infty} \left(-\frac{1}{\beta N}\, Z_N(\beta)\right)\, ,
$$
satisfies
\begin{equation}
\label{eq:conj}
\mathscr{F}\, \stackrel{?}{=}\, \max_{\mu \in \mathcal{M}(\Delta^{\infty})} \E[-\beta^{-1} \mathscr{S}(\xi) + \mathscr{E}(\xi)]\, .
\end{equation}
We would also hope that the arg-max is unique and is given by the limiting distribution of the random partition
structure 
$$
\left(\frac{X_{N}^{(1)}(\beta)}{Z_N(\beta)},\dots,\frac{X_{N}^{(2^n)}(\beta)}{Z_N(\beta)},0,0,\dots\right)\, ,
$$
where $X_{N}^{(1)}(\beta) \geq \dots \geq X_{N}^{(2^N)}(\beta)$ are the order statistics
$$
\sum_{n=1}^{2^N} \delta_{X_N^{(n)}(\beta)}\, =\, \sum_{n=1}^{2^N} \delta_{X_{N,n}(\beta)}\, .
$$
The reason for taking the maximum instead of the minimum is a deep insight of Giorgio Parisi which we will not discuss here.
(See for example \cite{MPV}.)

There are actually variational principles for the REM, but we do not know of this statement appearing before, although we also do not have a proof of this.
On the other hand, for the special case of Poisson-Dirichlet distributions $\operatorname{PD}(m,0)$ one can do the calculations
explicitly.
This follows because for any IID random variables $W_1,W_2,\dots \geq 0$, independent of $\xi \sim \operatorname{PPP}((0,\infty),mx^{-m-1}\, dx)$,
we have 
$$
\sum_{n=1}^{\infty} \delta_{W_n \xi_n}\, \stackrel{\mathcal{D}}{=}\, \sum_{n=1}^{\infty} \delta_{\E[W_n^m]^{1/m} \xi_n}\, .
$$
See, for example, \cite{RuzmaikinaAizenman} for an informative proof of this important stability property.
This implies that 
\begin{align*}
\E[\mathscr{E}(\xi)]\, &=\, -(\beta m)^{-1} \ln \E[e^{-\beta m Z_1}]\, =\, -\frac{\beta m}{2}\, ,\ \text{ and }\\
\E[\mathscr{S}(\xi)]\, &=\, m^{-1} \ln [B_1^m]\, =\, m^{-1} \ln 2\, .
\end{align*}
So restricting attention to the Poisson-Dirichlet distributons we would have
$$
\mathscr{F}(\beta)\, =\, \max_{m \in (0,1]} \left(-(\beta m)^{-1} \ln 2 - \beta m/2\right)\, .
$$
This gives the correct answer of course: for $\beta \leq \beta_c = \sqrt{2\ln 2}$, one takes $m=1$
and for $\beta>\beta_c$ one takes $m=\beta_c/\beta$ to obtain
$$
\mathscr{F}(\beta)\, =\, \begin{cases} -\beta^{-1} \ln 2 - \beta/2 & \text{ for $\beta\leq \beta_c$,}\\
-\beta_c & \text{ for $\beta > \beta_c$.}
\end{cases}
$$
Also, this leads to the correct partition structure $\operatorname{PD}(m(\beta),0)$.
Of course, we have heavily handicapped the problem since we already know that this is the solution
obtained by Derrida \cite{Derrida} and is even close to the original ideas.
On the other hand, it would be interesting to know if (\ref{eq:conj}) is true if one considers the full
class $\mathcal{M}(\Delta_{\infty})$.

\section*{Acknowledgments}

We are very grateful to Amir Dembo, who alerted us to the important paper of Shkolnikov.
We also gratefully acknowledge useful suggestions of Dmitry Panchenko and Paul Jung.
Part of this research was carried out during the workshop, “Statistical Mechanics
on Random Structures,” organized in Winter 2009 at the Banff International Research
Station (09w5055). We are grateful to the organizers and to BIRS for the accommodating
environment.

\end{document}